\documentclass[11pt,bezier,epsf]{article}
\usepackage{amsmath,amssymb,amsfonts,pstricks}
\usepackage{graphicx}
\newtheorem{theorem}{\bf Theorem}

\newtheorem{thm}{Theorem}[section]

\newtheorem{definition}[thm]{Definition}

\newtheorem{example}[thm]{Example}
\numberwithin{equation}{section}
\begin{document}
\vspace*{1cm}
\begin{center}
{{\bf \Large Index definitions for  nonlinear IAEs and DAEs: new classifications and numerical treatments}}

\vspace{1cm} {\bf B. Shiri,}
\\
{\it Faculty  of  Mathematical  Sciences,  University  of  Tabriz, Tabriz-Iran.\\
 E-mail:  shiri@tabrizu.ac.ir\\
Tel:  +98-0411-3391929 \\
Fax: +98-0411-3342102}\\
\end{center}
%\date{}%
%\dedicatory{}%
%\commby{}%
\vspace{1cm}
% ----------------------------------------------------------------
\hrule
\begin{abstract}
The  definition  of  index  for  differential  algebraic  equations  (DAEs)  or  integral  algebraic
equations  (IAEs)  in  the  linear  case  (time variable)  depends  only  on
the  coefficients  of  integrals  or   differential  operators   and  the
coefficients  of  the  unknown functions.
Is  this  possible  for  the  nonlinear  case?  In  this  paper  we  answer   this  question.
 In this paper, we  generalize  the index notion
for  the  nonlinear  case.  One  of  the  difficulties for  nonlinear  case,  is
its   dependence   on  the exact  solution  which  motivates  us  to  give  an  important  warning  to  whom  want  to
solve   DAEs  using  numerical  methods  such  as
Runge-Kutta{, multistep or}  collocation  methods.\\
 {\bf keywords:}{ Differential algebraic equation; Integral algebraic equation; Runge-Kutta methods; Collocation methods; {Multistep methods}.}\\
\end{abstract}
%%%%%%%%%%%%%%%%%%%%%%%%%%%%%%%%%%%%%%%%%%%%%%%%%%%%%%%%%%%%%%%%%%%%%%%%%%%%%%%%%%%%%%%%%%%%%%%%%%%%%%%%%%%%%%%%%%%%%%%
\section{Introduction}
We  consider  DAEs   and  IAEs  of  the  linear  forms
\begin{equation}\label{DAEs1}
A(t)y{'(t)}+B(t)y(s)ds=f(t), \ \ \ \ \ \ \ t\in I:=[0, T]
\end{equation}
and
\begin{equation}\label{IAEs1}
A(t)y(t)+\int_0^tk(t,s)y(s)ds=f(t), \ \ \ \ \ \ \ t\in I:=[0, T]
\end{equation}
with  $A, B \in  \mathbf{C}(I,{R}^{r\times r}),$
$f\in\mathbf{C}(I,{R}^{r}),$
$k\in \mathbf{C}(\mathbb{D},{R}^{r\times r}),$
 and  their  semi-nonlinear  forms
\begin{equation}\label{DAEsn}
A(t)y'(t)+F(t,y(t))=f(t), \ \ \ \ \ \ \ t\in I:=[0, T],
\end{equation}
and
\begin{equation}\label{IAEsn}
A(t)y(t)+\int_0^t\kappa(t,s,y(s))ds=f(t), \ \ \ \ \ \ \ t\in I:=[0,
T],
\end{equation}
where  $\kappa \in \mathbf{C}(\mathbb{D}\times{R}^{r},{R}^{r})$  and
$F \in \mathbf{C}(I\times{R}^{r},{R}^{r})$  with
$\mathbb{D}:=\{(t,s): 0\leq s\leq t\leq T\}$.
We  also  assume  that  $A(t)$  is  a  singular  matrix  with  constant  rank  for  all  $t\in I.$
Since  integrating  (\ref{DAEs1})  and  (\ref{DAEsn})
changes  them  to  IAEs  of  the  forms  (\ref{IAEs1})  and  (\ref{IAEsn}),
we  conclude  that  the  DAEs  inherit  many  properties  of  IAEs.
We  are   interested  in  translating  many  concepts  and  results  on  DAEs  to  the  poorly  studied  IAEs.

We  recall  the  DAEs  and  IAEs  classification   by  {the  index
notions.}  There  are  many  index  notions  by  considering  analytical  and
numerical  properties  of  solutions. The  differentiation  index  for
DAEs  goes  back  to  the  work  of  Campbell  \cite{Campbell,Brenan}  and
for  IAEs  to  the  work  of  Gear  \cite{Gear}.  Gear  introduced  the
differentiation  index  using  index  reduction  procedure  \cite{Gear}
to  obtain  the  existence  and  uniqueness  conditions  of  IAEs. These
definitions  of  index  are  extended  in  two  different  ways:  the
pioneering  works  of  Griepentrog(1991),  Reich(1991)  and  Rabier,
Rheinbold(1991)  \cite{Griepentrog,Reich,Rabier}  which  are  related  to
the  differential  geometric  concepts  of  regular  DAEs,  and  the
 work  of  Chistyakov  and  Bulatov  (\cite{Bulatov,Chistyakov})  (left index).  The
 global  index  which  is  an  extension  of  index  notion  for  linear  DAEs
 with  constant  coefficients  using  Kronecker  canonical  normal  form,  was  given  by  Gear  and  Petzold \cite{Petzold}.
In  \cite{Marz1},  Marz  has  shown  that  the  global  index  must  be  replaced  by  tractable  index.
{ There are two reasons in \cite{Marz1}; 1) there are not any way for relating this index to the nonlinear case, 2) unavailability of the necessary transforming matrices, except for some interesting cases.}
The  projection  type  index  definition  for  DAEs (tractable index)  which was
introduced  by  Marz  and  her  colleagues  \cite{Marz1},  can  be  extended
to  IAEs,  which  can  be  found  at   the  works  of  Brunner  and  Pishbin  \cite{herman,Hadizadeh}.
The  perturbation  index  defined  by  Hairer,  Lubich  and  Roche
\cite{Lubich}  has  important  role  in  analyzing   numerical  treatments
of  DAEs.  Finally,  strangeness  index  was  introduced  by  Kunkel  and  Mehrmann
\cite{Kunkel}.

This  paper  is  organized  as  follows:  in  section   $2,$  the  preliminary  definitions  of  index  for  the  linear  case  are  given,  also  a
theorem  for  existence  and  uniqueness  of  solution  of  IAEs  is  proved.  In  section   $3,$  after  defining index  for  nonlinear  case,  the  related  problems  are discussed.  Finally,  by  giving  numerical  experiments,
we  observe  the  cautions  that should  be  taken  for  solving  IAEs  and  DAEs.
%section $9,$ we illustrate the results of paper by giving some
%numerical experiments.

%%%%%%%%%%%%%%%%%%%%%%%%%%%%%%%%%%%%%%%%%%%%%%%%%%%%%%%%%%%%%%%%%%%%%%%%
\section{Index  definition  for  the  linear  case }
There are many definitions for index and one can refer to the reference in the introduction. Although, in this paper, we only use `rank degree'  index for linear case, but other index also can be used for generalization to the nonlinear case. To define `rank degree'  index, we need following preliminaries.
\begin{definition}\cite{Bulatov,Chistyakov}
 The  matrix  $A^{-}(t)$  is  called  a  semi-inverse  of  $A(t),$  if  it
 satisfies  the  equation
$$A(t)A^{-}(t)A(t) = A(t),$$
which  can   be  rewritten  as
$$V(t) A(t) = 0,$$
%where $E$ is $r\times r$ identity matrix
 with
\begin{equation}\label{s12}
V(t) = {E}-A(t)A^{-}(t),
\end{equation}
where  ${E}$  is  the  $r\times r$  identity  matrix.
 \end{definition}

 The  following  conditions  are  necessary  and  sufficient  for  the  existence  of  a  semi-inverse  matrix $A^{-}(t)$  with  elements  in $C^p({I},{R}^{r\times r})$ \cite{Chistyakov}:
\begin{enumerate}
  \item  the  elements  of $A(t)$   belong  to  $C^p({I},{R}^{r\times r});$
  \item  $\mbox{rank} A(t)=\mbox{constant},$ $\forall t\in {I}.$
\end{enumerate}
 \begin{definition}\cite{shiri}
 Suppose  that  $A\in C^{\nu}(I,R^{r\times r})$  and  $k\in C^{\nu}(D,R^{r\times r}).$  Let
 $$A_0\equiv A, \ \ \ \ k_0\equiv k,$$
 $$\Lambda_i y=\frac{d}{dt}\left(({E}-A_i(t)A^-_i(t))y\right)+y,$$
 $$A_{i+1}\equiv A_{i}+({E}-A_i(t)A^-_i(t))k_i(t,t), \ \mbox{and} \ k_{i+1}=\Lambda_i k_i.$$
  Then,  we  say   the  `rank degree'
index  of  $(A, \ k)$  is  $\nu,$  if
$$\mbox{rank} A_i(t)=\mbox{constant}, \ \ \ \forall t\in{I}, \ \ \mbox{for} \ \ i=0 ,\ldots, \nu,$$
$$\det A_i=0,\ \ \mbox{for} \ \ i=0 ,\ldots, \nu-1, \ \  \det A_{\nu}\neq 0.$$
Moreover,  we  say   the  `rank degree'  index  of {the} system
(\ref{IAEs1})  is  $\nu$ ($ind_r=\nu$),
 if  in  addition  to  the  above  conditions,  we  have  $f\in C^{\nu}(I,R^r)$  {and hence, we can define}
$$F_{i+1}\equiv \Lambda_i F_i, \ {\mbox{with}} \ F_0\equiv f, \ i=0, ..., \nu-1.$$
 \end{definition}
 Now, we can state the following theorem of  the uniqueness and existence for higher index IAEs.
 \begin{theorem}\label{sth20}\cite{shiri}
 Suppose  the following conditions are satisfied for (\ref{IAEs1}):
\begin{enumerate}\label{existence2}
          \item  $ind_r=\nu\geq 1,$
          \item $A_i(t)\in \mathbf{C}^1({I},{R}^{r\times r}),$
           $F_i(t)\in \mathbf{C}^1 ({I},{R}^{r})$ and $k_i\in \mathbf{C}^1(\mathbb{D},{R}^{r\times r})$  for $i=1, \cdots, \nu,$
          \item      $A_i(0)A_{\nu}^{-1}(0)F_{\nu}(0)=F_i(0)$ for $i=0, \cdots \nu-1$ (consistency conditions).
        \end{enumerate}
        Then  the system (\ref{IAEs1}) has a unique solution on $I.$
 \end{theorem}

  \begin{definition}\label{globald}\cite{Petzold}
We say the  `rank degree'  index  of the system (\ref{DAEs1})  is   $\nu,$ if  the  `rank degree'  index  of  it's corresponding IAE
\begin{equation}\label{DAE2IAE}
A(t)x(t)+\int_0^t(B(s)-A'(s))x(s)ds=\int_0^tq(s)ds.
\end{equation}
be  $\nu.$
\end{definition}

\section{index definition for the nonlinear case}
Among  the  existing  definitions  for  index,  there  are  definitions
which  are  independent  of  the  linear  case:
\begin{definition} \cite{Campbell,Brenan}
We say the differentiation index of the system (\ref{DAEsn}) is $\nu$ ($ind_d=\nu$), if $\nu$ is
the minimum possible number of derivatives (\ref{DAEsn}) to obtain a system of the
ordinary differential equations(ODE).
\end{definition}
\begin{definition}\cite{Gear}
We say the differentiation index of the system (\ref{IAEsn}) is $\nu$ ($ind_d=\nu$), if $\nu$ is the
minimum possible number of derivatives (\ref{IAEsn}) to obtain
a system of the second kind Volterra integral equations(SVIE).
\end{definition}
\begin{definition}\label{perturb}\cite{Lubich}
The  equation $F(Y',Y)=0$  has  perturbation  index  $\nu$  ($ind_p=\nu$)
along  a  solution  $Y$  on  $I,$  if  $\nu$  is  the  smallest  integer  such
that  for  all  functions  $\widehat{Y}$  having  the residual
$$F(\widehat{Y}',\widehat{Y})=\delta(x),$$
there  exists   an  estimate
$$\|\widehat{Y}(x)-Y(x)\|\leq C(\|\widehat{Y}(0)-Y(0)\|+\max_{x\in I}\|\delta(x)\|+\ldots+\max_{x\in I}\|\delta^{(\nu-1)}(x)\|)$$
on  $I,$  whenever  the  expression  on  the  right-hand  side  is
sufficiently  small.  Here,  $C$  denotes  a  constant  which  depends
on  $F$  and  the  length  of  the  interval.
\end{definition}

The  perturbation  index  defined  by  Hairer,  Lubich  and  Roche
\cite{Lubich}  has  an  important  role  in  analyzing  numerical  treatment
of  DAEs,  where  the  inequality
$$\|\widehat{Y}(x)-Y(x)\|\leq C(\max_{x\in I}\|\delta(x)\|+\ldots+\max_{x\in I}\|\delta^{(\nu)}(x)\|)$$
should  be  considered  in  above  definition  for  IAEs.

Now  the  question  is:  how   can  we   extend  the   index  definitions   of
previous  sections  to  the  nonlinear  case?

Suppose  $u$  denotes   an  approximate  solution  of  the  nonlinear
equation  (\ref{DAEsn})  or  (\ref{IAEsn}),  and  define  the  residuals
\begin{equation}\label{DAEsnd}
R_1(t):=A(t)u'(t)+F(t,u(t))-f(t)
\end{equation}
or
\begin{equation}\label{IAEsnd}
R_2(t):=A(t)u(t)+\int_0^t\kappa(t,s,u(s))ds-f(t),
\end{equation}
which  are  supposed  to  be  very  small.  Subtracting  (\ref{DAEsn})  or
(\ref{IAEsn})  from  (\ref{DAEsnd})  or (\ref{IAEsnd}),  we  obtain
\begin{equation}\label{DAEsnd}
R_1(t)=A(t)(u'(t)-y'(t))+F(t,u(t))-F(t,y(t))
\end{equation}
or
\begin{equation}\label{IAEsnd}
R_2(t)=A(t)(u(t)-y(t))+\int_0^t\kappa(t,s,u(s))-\kappa(t,s,y(s))ds.
\end{equation}
Letting  $e(t):=u(t)-y(t),$   and  supposing  that  $F$  and  $\kappa$
have  continuous  derivatives   with  respect  to  the  second  and  third
variables  respectively,  we  have
\begin{equation}\label{DAEsndl}
R_1(t)=A(t)e'(t)+F_y(t,\eta(t))e(t)
\end{equation}
or
\begin{equation}\label{IAEsndl}
R_2(t)=A(t)e(t)+\int_0^t\kappa_y(t,s,\eta(s))e(s)ds,
\end{equation}
where  $\eta(s)$  is  an  appropriate  function  arises  from  applying  the
mean  value  theorem  which  is  in  a  neighborhood  of  the  exact  solution  $y$.
The  equations  (\ref{DAEsndl})  and  (\ref{IAEsndl})  are  of  the  linear
form  (\ref{DAEs1})  and  (\ref{IAEs1}).

Now,  the  strategy  for  defining  an  index  for  the  nonlinear  case
(\ref{DAEsn})  or  (\ref{IAEsn})  is  the  same  as  in  linear  cases
(\ref{DAEsndl})  or  (\ref{IAEsndl}).

\begin{definition}
We  say  that  the  index  for  the  nonlinear  equations  (\ref{IAEsn})  or
(\ref{DAEsn})  is  $\nu,$   if  there  exists  a  neighborhood  of  the  exact
solution  $y$  ($N_\epsilon(y)$),  in  which  the  index  of  the  linear
equations  (\ref{IAEsndl})  or  (\ref{DAEsndl})  is  $\nu.$
\end{definition}

This  dependency  of  index  definition  to  the  exact  solution  (that is unknown)   is  the  difficulty  of  this  definition  and  definition
\ref{perturb}.  We  will  show  by  examples,   how  we  can   use  this
definition  effectively.

We  divide  the  nonlinear  IAEs  and  DAEs  to  the  following  classes:
\begin{description}
  \item[well structure:]  the  IAEs  or  DAEs  that  their  indices  do  not %
  depend  on  the  unknown  variables  (component  of  $y$).
  \item[free structure:]  the  IAEs  or  DAEs  that  their  indices
   depend  on  the  unknown  variables.   %
\end{description}
Fortunately,  most  of  IAEs   and  DAEs  in  application  are  well
structure.   For   example   the  DAEs  in  \cite{Lubich}  are  well
structure  and  for  every  neighborhood  of  $y$  their  indices  remain
constant. The IAEs  of  Hessenberg  type
 \begin{equation}\label{ind2qn}
 \begin{split}
&\left[
  \begin{array}{cccc}
  A_{1,1}(t)& \ldots & A_{1,\nu-1}(t) & 0\\
  \vdots& \ddots& \vdots & \vdots\\
 A_{\nu-1,1}(t)& \ldots &0&0\\
 0& \ldots & 0& 0
  \end{array}
\right]\left[\begin{array}{c}
                             y_1(t) \\
                             y_2(t) \\
                             \vdots \\
                             y_{\nu}(t)
                           \end{array}
                         \right]+\\
                         &\int_0^{t}\left[
                    \begin{array}{c}
                      k_{1}(t,s,y_1(s),y_2(s),\ldots,y_{\nu}(s)) \\
                        \vdots\\
                      k_{\nu-1}(t,s,y_1(s),y_2(s))\\
                      k_{\nu}(t,s,y_1(s))
                    \end{array}
                  \right]ds
                  =\left[
                            \begin{array}{c}
                              f_1(t) \\
                              f_{2}(t) \\
                              \vdots \\
                             f_{\nu}(t)
                            \end{array}
                          \right],
                          \end{split}
\end{equation}
are  well  structure  with   index   $\nu$,
 if  $\prod_{i=1}^{\nu}k_{i,y_{\nu+1-i}}(t,t,y_1(t),...,y_{\nu+1-i}(t))$  is  invertible  with  a  bounded  inverse  on I.
\begin{example}\label{ex31}
The  DAE   of  the  form
\begin{equation}
\begin{split}
y'&=y^2+e^{y}+z, \\
0&=e^{y}+\sin(t),
\end{split}
\end{equation}
is  a  well  structure  DAE  of   index  2,  since
 $$det(F_y)= \left[
  \begin{array}{cc}
    2y+e^{y} & 1 \\
   e^{y}& 0\\
  \end{array}
\right]=e^{y}>0, \ \ \mbox{for all } y \mbox{ and } z.
$$
\end{example}

To  speak  about  the  free  structure  IAEs  or  DAEs,  we   divide  them
again  to  the  following  classes:
\begin{description}
  \item[the dependent form:]  the  IAEs  or  DAEs  that  their  indices
  depend  on  the  exact  solution  on  the  given  interval.
  \item[the independent form:]  the IAEs  or  DAEs  that  their  indices  do  not
  depend  on  the  exact  solution  on  the  given  interval.
\end{description}
A  condition  is  called  critical,  if  it  changes
the  index  in  a  free  structure  DAE  or  IAE.  Also,  a  point  is  called   critical  if  at  which  the  exact
solution  satisfies  the  critical  conditions  of  the dependent
form  DAE  or  IAE.

\begin{example}\label{ex32}
 Let
$$A(t)=\left(
         \begin{array}{cc}
           1 & 0 \\
           0 & 0 \\
         \end{array}
       \right),
$$
$$F(t,y)=\left(
           \begin{array}{c}
             -y_1^2-e^{y_2} \\
             -y_1y_2 \\
           \end{array}
         \right),
$$
and
$$f(t)=\left(
           \begin{array}{c}
             \cos^2(t)-e^t-sin(t) \\
             -t\cos(t) \\
           \end{array}
         \right).
$$
 Then a DAE of the form  (\ref{DAEsn}) has the exact  solution
  $$
  y(t)=\left(
         \begin{array}{c}
           \cos(t) \\
           t \\
         \end{array}
       \right)
  $$
 and we have
  $$F_y(t,\eta)=\left[
          \begin{array}{cc}
            -2\eta_1 & -e^{\eta_2} \\
            -\eta_2 & -\eta_1 \\
          \end{array}
        \right].
$$
This  is  a  free  structure  DAE,  and  it  is  obvious  that  for  $\eta_1\ne
0$  the  index  is  1,  and  for  $\eta_2\ne 0, \eta_1= 0$  the  index  is
$2.$  Now  suppose  the  interval  of  solution  is  $I_1=[.5, 1].$  In  this
interval,  the  exact  solution  doesn't  vanish  and  we  can  find
$\epsilon>0$  such  that  $[\eta_1,\eta_2]^T\in
N_{\epsilon}([\cos(t),t]^T)$  and  $\eta_1\ne 0$  for  all  $t\in I_1.$
Hence,  on  $I_1$  the  index  is  do not change  and  the  DAE  is in the  independent
form.  But  for  the  interval  $I_2=[1,\ 2],$  we
have  $y_1(\pi/2)=0$   and  the  index  changes  from  1  to  2  at  $t=\pi/2$.  Hence,
this  DAE  has  a  dependent  form on  $I_2$.
\end{example}
\begin{example}\label{ex33}
Let  $A(t)$  and  $F(t,y)$  be  the  same  matrices  introduced  in  the  previous
example,  and  let
$$f(t)=\left(
           \begin{array}{c}
             -e^{2t} \\
              -te^{t}\\
           \end{array}
         \right).
$$
Then  the  exact  solution  is  $[e^{t},t]^T.$  Obviously,  the  index  of  this  system  remain
constant  on  $I_1$  and
$I_2$  and  on  every other  bounded  interval  and  it  is  equal  to  1.  This  free  structure  DAE  has  an
independent  form   of  index  1.
\end{example}
\begin{example}\label{ex34}
Consider  an  IAE  of  the  form  (\ref{IAEsn}).  Let
$$A(t)=\left(
         \begin{array}{cc}
           1 & 0 \\
           0 & 0 \\
         \end{array}
       \right)
$$
and
$$\kappa(t,s,y(s))=\left(\begin{array}{c}
                     (y_1^2+2)y_2+e^{y_2} \\
                     y_1^2
                   \end{array}\right)
$$
and  let  $f$   be  a  function  such  that  the  exact  solution  to be
$y(t)=[e^{t},t]^T.$  For  this  IAE,  we  have
  $$\kappa_y(t,s,\eta(s))=\left[
          \begin{array}{cc}
            2\eta_1\eta_2 & (\eta_1^2+2)+e^{\eta_2} \\
            2\eta_1 &0 \\
          \end{array}
        \right].
$$
Then the  index  of  its  corresponding  linear  IAE  is   2,  if  $\eta_1\neq
0,$  and it  is   $0,$  if  $\eta_1=0.$  Hence,  this  IAE  has  a  free
structure.  It  is  obvious  that  this  IAE  has  an  independent  form  of
index  2.
\end{example}

\begin{example}\label{ex35}
Let  $A(t)$  and  $\kappa(t,s,y(s))$  be  the  same  matrix  introduced  in
Example  \ref{ex33}  and  let
$$f(t)=\left(
           \begin{array}{c}
             \cos(t) - \frac{\sin(2)}{4} - e + e^{t} + \frac{t\sin(2t)}{4} - \frac{\sin(t)^2}{4} + \frac{\sin(1)^2}{4} + \frac{5t^2}{4} - \frac{5}{4} \\
              \frac{t}{2} + \frac{\sin(2t)}{4} - \frac{\sin(2)}{4} - \frac{1}{2}\\
           \end{array}
         \right).
$$
Then  the  exact  solution  is  $[\cos(t),t]^T.$   This  IAE  has  a free
structure  dependent  form  with  index  changing  at  the  point  $t=\pi/2.$
\end{example}

\section{Numerical experiments}

\begin{center}
\begin{figure}
  \includegraphics[width=.8\textwidth]{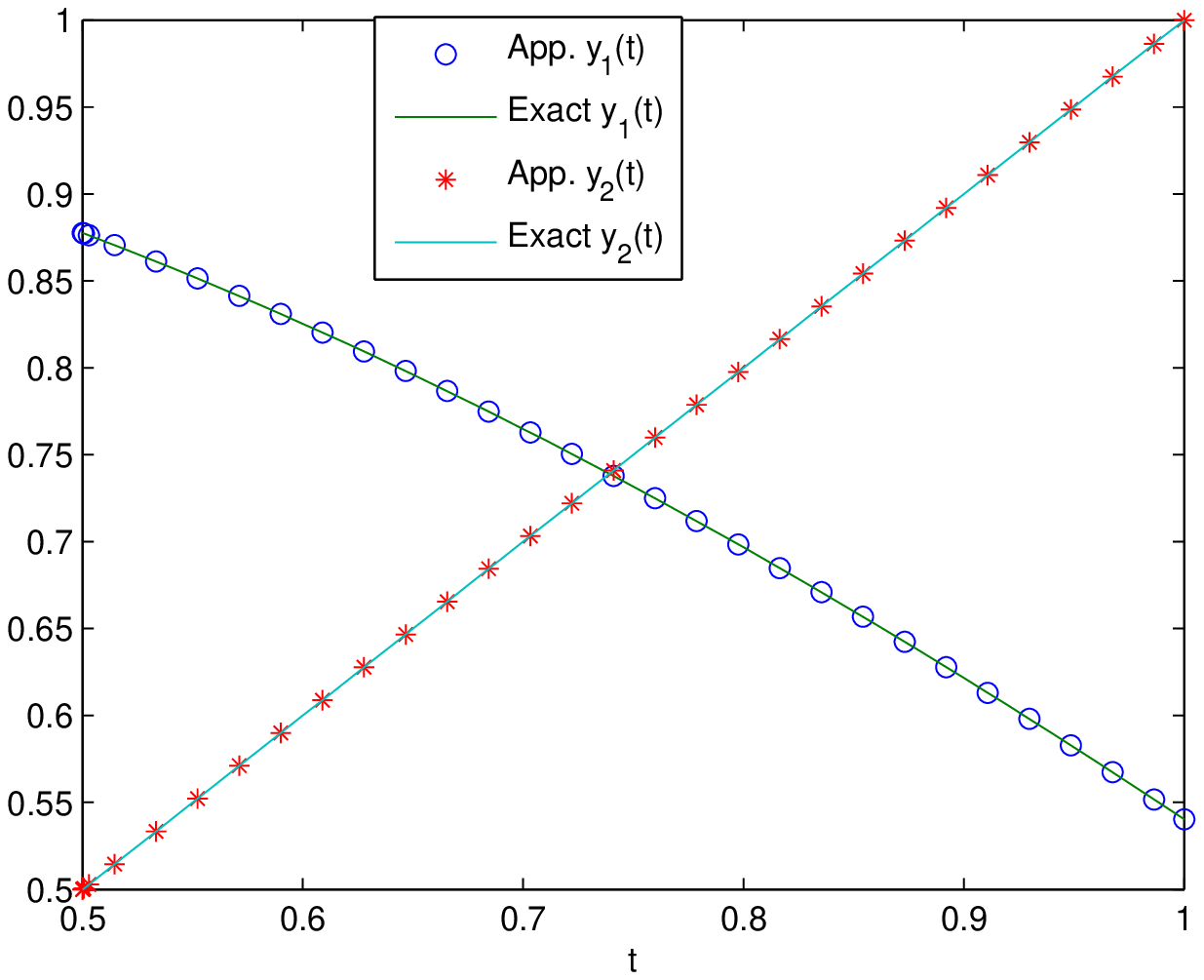}
 \caption{The  numerical  and  exact  solutions  of  the  free  structure
independent  form   DAE   in Example  \ref{ex36}.}
 \label{figg1}
\end{figure}
\end{center}
\begin{center}
\begin{figure}
  \includegraphics[width=.75\textwidth]{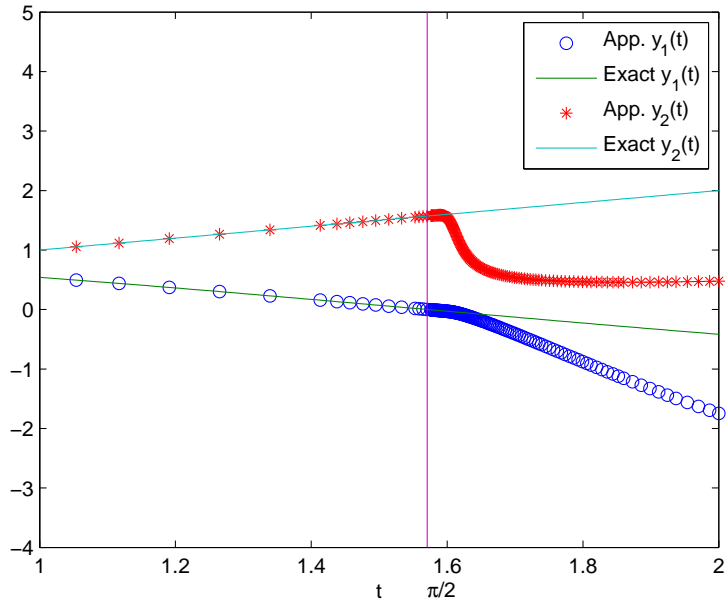}
 \caption{The  numerical  and  exact  solutions  of  the
free  structure  dependent  form   DAE  in  Example \ref{ex36}.  This  figure
shows  that  after  the  point  $t=\pi/2$  the  numerical  solution  does  not
converge  to  the  exact  solution. This results are obtained using the command  `ode15s' with the options $`RelTol'=1e-6$ and $`AbsTol'=[1e-8, 1e-8].$}
 \label{figg2}
\end{figure}
\end{center}
\begin{center}
\begin{figure}
  \includegraphics[width=.75\textwidth]{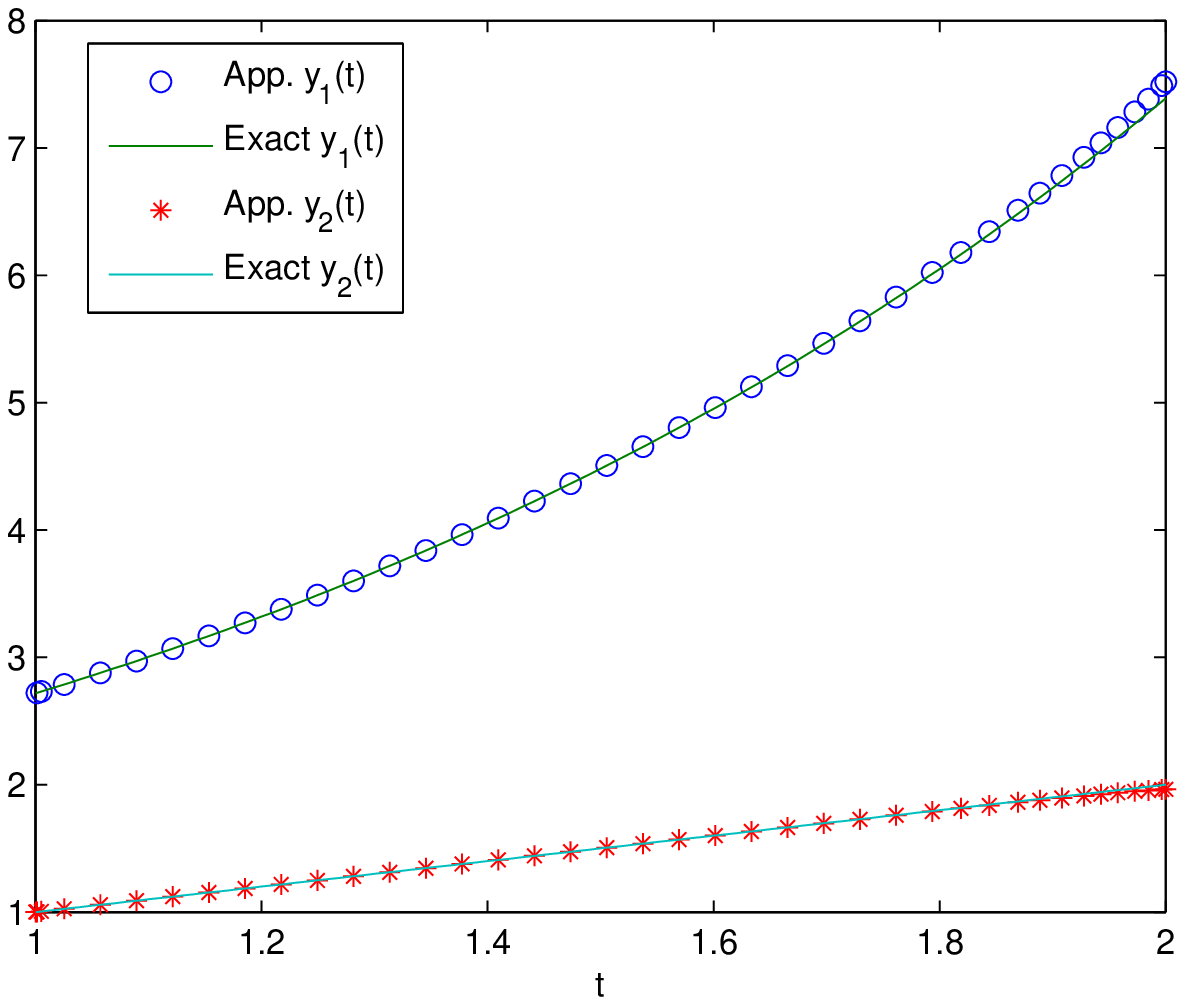}
  \begin{center}
 \caption{The  numerical  and  exact  solutions  of  free  structure
independent  form   DAE  in  Example  \ref{ex36}.}
\end{center}
 \label{figg3}
\end{figure}
\end{center}
\begin{center}
\begin{figure}
  \includegraphics[width=.82\textwidth]{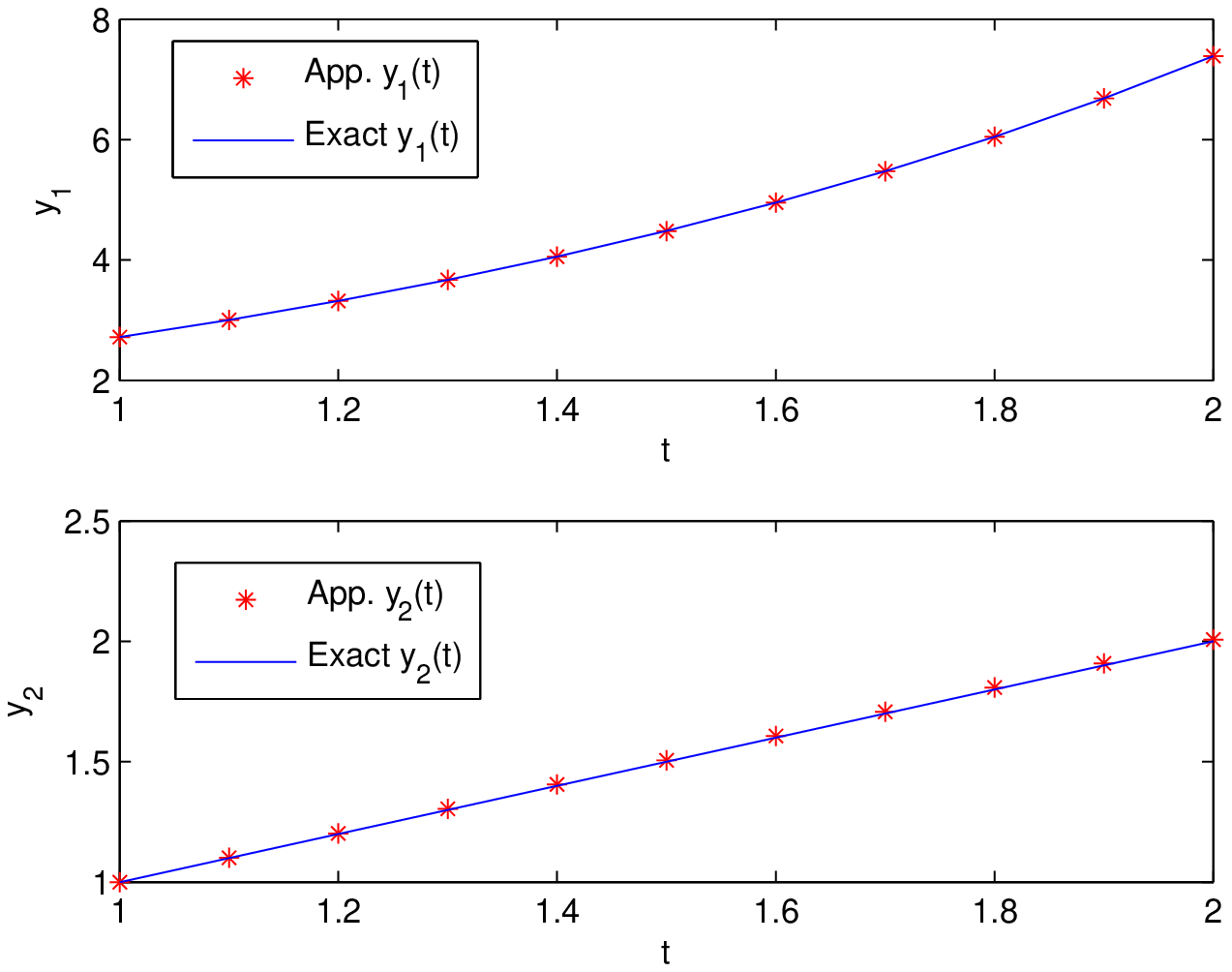}
 \caption{The  numerical  and  exact  solutions  of  free  structure
independent  form   IAEs  in  Example   \ref{ex37}.}
 \label{figg4}
\end{figure}
\end{center}
\begin{center}
\begin{figure}
  \includegraphics[width=.82\textwidth]{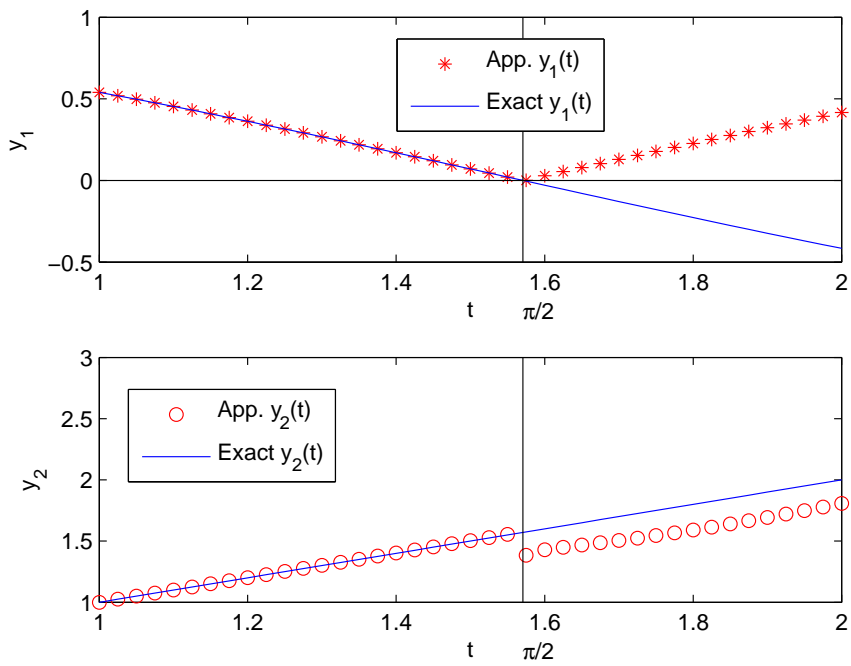}
  \caption{The  numerical  and  exact  solutions  of
independent  form  free  structure  IAEs  in  Example  \ref{ex37}.  We  can
observe  that  numerical  solution  doesn't  converge  to  the  exact  solution
after  the  point  $\pi/2$.}
 \label{figg5}
\end{figure}
\end{center}

 We  use  the  MATLAB  codes  for  solving  DAEs  and  the
collocation  methods  introduced  in  \cite{herman}  on  the  piecewise
polynomial  space  for  IAEs.

\begin{example}\label{ex36}
We  solve   the  free  structure   DAE  of  Example  \ref{ex32}  on  $I_1$  and
$I_2$  using  the  command  `ode15s'  in MATLAB.  The  index  on  $I_1$
remains  constant  and  we  see  that  the  numerical  solution  gives  a  good
approximation  of  the  exact  solution  (see  Figure  \ref{figg1}).  But
when  we  use  `ode15s'  for  this  DAE   on   $I_2,$  dependent  form,
the  behaviors  of  the  exact   and  approximate  solutions  change
after  $t=\frac{\pi}{2}$  and  this  is  because  of  changing  the  index  at
this  point  (see  Figure  \ref{figg2}).  This  example  shows  that  we
should  make  caution,  when  we  solve  a  free  structure  DAE.  Figure
\ref{figg3},  shows  the  numerical  and  exact  solutions  of  the  free
structure  independent  form   DAE  in  Example  \ref{ex33}.  This  figure  shows  that  the   MATLAB  code  `ode15s'
can  be  used  for  solving  free  structure   DAEs,  although  some  care  is  needed.

\end{example}
\begin{example}\label{ex37}
We  solve  the  free  structure  IAEs  of  Examples  \ref{ex34}  and
\ref{ex35}  using  a  continuous  piecewise  collocation  method  on  the
space  $\mathcal{S}_{m}^{(0)}$  with $c=[0,.7,.9],$ and $h=0.025$ (see \cite{herman}).  The  same  phenomenon
of  the  previous  Example  is  observed  again. Figures   \ref{figg4}  and
\ref{figg5},  show  the  numerical  and  exact  solutions  of these  IAEs   on  $[1, 2].$
\end{example}

In  practice,  we  do  not  have  the  exact  solution  to  find  the  points
that  make  a  dependent  form  IAE  or  DAE.  Thus  the  problem  is  that:  how
 can  one  determine  the  interval  of  convergence?  The answer
 is  that,
we  have  the  data  of  numerical  solution,  using  these  data  and
critical  conditions  one  can  find  the  critical  points  and  end  the  solver
or  be  careful  about  the  approximate  solution  after  the  first
critical  point.  In  Example  \ref{ex32}  the  condition  $y_1=0$  is  a
critical  condition.  We  can  see  from  Figure  \ref{figg2}  that  the  numerical
values  of  $y_1$  tends  to  zero  as  t  tends  to  $\pi/2.$  So,  we
 caution  about  the  solution  after  this  point.  The  same  thing
are  observed  in  the  numerical  solution  of  Example  \ref{ex34}  in
figure  \ref{figg5}.

\section{Conclusion}
The  difficulty  of  index  definitions  in   nonlinear  DAEs  and  IAEs  lead
to  dividing  the  IAEs  and  DAEs  to  three  classes:   well  structure,
free  structure  independent  form  and  free  structure  dependent  form.
We  observed  that  the  dependency  of  index  definition  to
unknown  solution  can  not  be  removed,  since  it  affects  on  the  numerical
treatments  of  solutions.  We  concluded  that  we  can  use  the data of numerical
   solutions, but we must be careful  in  solving  nonlinear  DAEs  and  IAEs.

\bibliographystyle{amsplain}
%\bibliography{}
\end{document}